\documentclass[12pt,epsfig]{article}
\usepackage{graphicx}
\usepackage{array}
\usepackage{comment}
\usepackage{amsmath,amssymb}
\usepackage[latin1]{inputenc}
\begin{document}
\title{Algebras with ternary law of composition combining $Z_2$ and $Z_3$ gradings} %%%%% Ispravil Z_3 na Z_2

\author{V. Abramov\footnote{Institute of Mathematics,
University of Tartu, Liivi 2, Tartu 50409, Estonia,
E-mail: viktor.abramov@ut.ee, olga.liivapuu@gmail.com}, R. Kerner\footnote{LPTMC, Tour 23-13, 5-e \'etage, Bo\^{i}te 121, 
4, Place Jussieu, 75252 Paris Cedex 05, France,~~E-mail: richard.kerner@upmc.fr} and O. Liivapuu$^*$}

\date{ }

\maketitle

{\large {\bf Summary}}

\vskip 0.3cm
\indent
In the present article we investigate the possibility of combining the usual Grassmann algebras with
their ternary $Z_3$-graded counterpart, thus creating a more general algebra with coexisting quadratic and cubic
constitutive relations.

We recall the classificaton of ternary and cubic algebras according to the symmetry properties of ternary
poducts under the action of the $S_3$ permutation group. Instead of only two kinds of binary algebras,
symmetric or antisymmetric, here we get {\it four} different generalizations of each of these two cases.

Then we study a particular case of algebras generated by two types of variables, $\xi^{\alpha}$ and $\theta^A$,
satisfying quadratic and cubic relations respectively, $\xi^{\alpha} \xi^{\beta} = - \xi^{\beta} \xi^{\alpha}$
and $\theta^A \theta^B \theta^C = j \theta^B \theta^C \theta^A$, $j = e^{\frac{2 \pi i}{3}}$.

The invariance group of the generalized algebra is introduced and investigated.

\newpage

\section{Classification of ternary and cubic algebras}

In \cite{Kerner1991}, \cite{Kerner1992}, \cite{Abramov1997}, \cite{RKOS} certain types of $Z_3$-graded ternary and cubic algebras
have been introduced and investigated. In \cite{VainKer} a general classification of $n$-ary algebras was given.

The usual definition of an algebra involves a linear space ${\cal{A}}$ (over real or complex numbers) endowed with a {\it binary} constitutive relations
\begin{equation}
{\cal{A}} \times {\cal{A}} \rightarrow {\cal{A}}.
\label{algebradef}
\end{equation}
In a finite dimensional case, dim ${\cal{A}}$ = N, in a chosen basis ${\bf e}_1, {\bf e}_2, ..., {\bf e}_N$, the constitutive relations (\ref{algebradef})
can be encoded in {\it structure constants} $f^k_{ij}$ as follows:
\begin{equation}
{\bf e}_i  {\bf e}_j = f^k_{ij} \, {\bf e}_k.
\label{structure2}
\end{equation}
With the help of these structure constants all essential properties of a given algebra can be expressed, e.g. they will define a {\it Lie algebra} if they
are antisymmetric and satisfy the Jacobi identity:
\begin{equation}
f^k_{ij} = - f^k_{ji}, \; \; \; \; f^k_{im} f^m_{jl} + f^k_{jm} f^m_{li} + f^k_{lm} f^m_{ij} = 0,
\label{structureLie}
\end{equation}
whereas an abelian algebra will have its structure constants symmetric, $f^k_{ij} = f^k_{ji}$.

Usually, when we speak of algebras, we mean {\it binary algebras}, understanding that they are defined via {\it quadratic}
constitutive relations (\ref{structure2}). On such algebras the notion of $Z_2-grading$ can be naturally introduced.
An algebra ${\cal{A}}$ is called a $Z_2-graded$  algebra if it is a direct sum of two parts, with symmetric (abelian) and
anti-symmetric product respectively,
\begin{equation}
{\cal{A}} = {\cal{A}}_0 \oplus {\cal{A}}_1,
\label{Agraded}
\end{equation}
with {\it grade} of an element being $0$ if it belongs to ${\cal{A}}_0$, and $1$ if it belongs to ${\cal{A}}_1$.
Under the multiplication in a $Z_2$-graded algebra the grades add up reproducing the composition law of the $Z_2$
permutation group: if the grade of an element $A$ is $a$, and that of the element $B$ is $b$, then the grade of their
product will be $a+b$ {\it modulo 2}:
\begin{equation}
{\rm grade} (AB) = {\rm grade} (A) + {\rm grade} (B).
\label{ABBAgraded}
\end{equation}
A $\mathbb Z_2$-graded algebra is called a $\mathbb Z_2$-graded commutative if for any two homogeneous elements $A,B$ we have
\begin{equation}
AB = (-1)^{ab} BA.
\label{graded commutative}
\end{equation}
It is worthwhile to notice at this point that the above relationship can be written in an alternative form, with all
the expressions on the left side as follows:
\begin{equation}
A B - (-1)^{a\,b} B A = 0, \; \; {\rm or} \; \; A B + (-1)^{(a\,b + 1)} B A = 0
\label{comanticom5}
\end{equation}
The equivalence between these two alternative definitions of commutation (anticommutation) relations inside a $Z_2$-graded algebra
is no more possible if by analogy we want to impose {\it cubic} relations on algebras with $Z_3$-symmetry properties, in which
the non-trivial cubic root of unity, $j = e^{\frac{2 \pi i}{3}}$ plays the role similar to that of $-1$ in the binary relations
displaying a $Z_2$-symmetry.

The $Z_3$ cyclic group is an abelian subgroup of the $S_3$ symmetry group of permutations of three objects. The $S_3$ group contains {\it six } elements,
including the group unit $e$ (the identity permutation, leaving all objects in place: $(abc) \rightarrow (abc)$), the two cyclic permutations
$$ (abc) \rightarrow (bca) \; \; \; {\rm and} \; \; \; (abc) \rightarrow (cab),$$
and three odd permutations,
$$(abc) \rightarrow (cba), \; \; \; (abc) \rightarrow (bac) \; \; \; {\rm and} \; \; \; (abc) \rightarrow (acb).$$

There was a unique definition of {\it commutative} binary algebras given in two equivalent forms,
\begin{equation}
xy + (-1) yx = 0 \; \; \; \; {\rm or} \; \; \; \; \; xy = yx.
\label{comm2}
\end{equation}
In the case of cubic algebras \cite{VainKer} we have the following four generalizations of the notion of {\it commutative} algebras:

\vskip 0.2cm
\begin{itemize}
\item[a)]
Generalizing the first form of the commutativity relation (\ref{comm2}), which amounts to replacing the $-1$ generator of $Z_2$ by
$j$-generator of $Z_3$ and binary products by products of three elements, we get
\begin{equation}
S: \; \; x^{\mu} x^{\nu} x^{\lambda} + j \;  x^{\nu} x^{\lambda} x^{\mu} + j^2 \; x^{\lambda} x^{\mu} x^{\nu} = 0,
\label{defS}
\end{equation}
where $j = e^{\frac{2 \pi i}{3}}$ is a primitive third root of unity.
\item[b)]
Another primitive third root, $j^2 = e^{\frac{4 \pi i}{3}}$ can be used in place of the former one; this will define thr conjugate algebra
${\bar{S}}$, satisfying the following cubic constitutive relations:
\begin{equation}
{\bar S} : \; \; x^{\mu} x^{\nu} x^{\lambda} + j^2 \;  x^{\nu} x^{\lambda} x^{\mu} + j \; x^{\lambda} x^{\mu} x^{\nu} = 0.
\label{defSbar}
\end{equation}
Clearly enough, both algebras are infinitely-dimensional and have the same structure. Each of them is a possible generalization of infinitely-dimensional
algebra of usual commuting variables with a finite number of generators. In the usual $Z_2$-graded case such algebras are just polynomials in variables
$x^1, \; x^2, ... x^N$; the algebras $S$ and ${\bar{S}}$ defined above are also spanned by polynomials, but with different symmetry properties, and
 as a consequence, with different dimensions corresponding to a given power.
\item[c)]
Then we can impose the following ``weak" commmutation, valid only for cyclic permutations of factors:
\begin{equation}
S_1 : \; \; x^{\mu} x^{\nu} x^{\lambda} = x^{\nu} x^{\lambda} x^{\mu} \neq x^{\nu} x^{\mu} x^{\lambda},
\end{equation}
%For a finite algebra of this type with $N$ generators the dimensions of the subsequent subspaces of polynomials of degree $k$ can be defined
% via the following {\it generating function} of one real vaiable $t$:
%\begin{equation}
%F(t) = \left( \frac{1}{1-t} \right)^N +(1+t)^N -1 - N = {\displaystyle{\sum_{k=0}^{\infty}}} f_k t^k = 1 + Nt + N^2 t^2 + \frac{2N (N^2+2)}{6} \, t^3 +...
%\end{equation}
\item[d)]
Finally, we can impose the following ``strong" commutation, valid for arbitrary (even or odd) permutations of three factors:
\begin{equation}
S_0: \; \; \; x^{\mu} x^{\nu} x^{\lambda} = x^{\nu} x^{\lambda} x^{\mu} = x^{\nu} x^{\mu} x^{\lambda}
\end{equation}
\end{itemize}
%The dimensions of the consecutive sets of polynomials of degree $p$ are given by the following generating function:
% \begin{equation}
%G(t) = \left( \frac{1}{1-t} \right)^N + \frac{N (N-1)}{2} \, t^2 = {\displaystyle{\sum_{k=0}^{\infty}}} g_k t^k = 1 + NT + N^2 t^2 + \frac{N (N+1)(N+2)}{6} \, t^3 +...
%\end{equation}
The four different associative algebras with cubic commutation relations can be represented in the following diagram, in which all arrows correspond
to {\it surjective homomorphisms}. The commuting generators can be given the common grade $0$.
\vskip 0.2cm
\centerline{$ S \hskip 1.5cm {\bar{S}}$}
\centerline{$ \searrow \hskip 0.5cm \swarrow$}
\centerline{$S_1$}
\centerline{$\downarrow$}
\centerline{$S_0$}
\vskip 0.2cm
Let us turn now to the $Z_3$ generalization of anti-commuting generators, which in the usual homogeneous case with $Z_2$-grading define Grassmann algebras.
Here, too, we have four different choices:
\vskip 0.2cm
\begin{itemize}
\item[a)] The ``strong" cubic anti-commutation,
\begin{equation}
\Lambda_0 : \; \; \; {\displaystyle{\Sigma_{\pi \in S_3} \; \theta^{\pi (A)} \theta^{\pi (B)} \theta^{\pi (C)}}} = 0,
\end{equation}
i.e. the sum of {\it all} permutations of three factors, even and odd ones, must vanish.
\item[b)]
The somewhat weaker ``cyclic" anti-commutation relation,
\begin{equation}
\Lambda_1 : \; \; \; \theta^A \theta^B \theta^C + \theta^B \theta^C \theta^A + \theta^C \theta^A \theta^B = 0,
\end{equation}
i.e. the sum of {\it cyclic} permutations of three elements must vanish. The same independent relation for the odd combination $\theta^C \theta^B \theta^A$
holds separately.
\item[c)]
The $j$-skew-symmetric algebra:
\begin{equation}
\Lambda : \theta^A \theta^B \theta^C = j \; \theta^B \theta^C \theta^A.
\end{equation}
and its conjugate algebra ${\bar{\Lambda}}$, isomorphic with $\Lambda$, which we distinguish by putting a bar on the generators and using dotted indices:
\item[d)]
The $j^2$-skew-symmetric algebra:
\begin{equation}
{\bar{\Lambda}} : \; \; \; {\bar{\theta}}^{\dot{A}} {\bar{\theta}}^{\dot{B}} {\bar{\theta}}^{\dot{C}} =
j^2 {\bar{\theta}}^{\dot{B}} {\bar{\theta}}^{\dot{C}} {\bar{\theta}}^{\dot{A}}
\end{equation}
Both these algebras are finite dimensional. For $j$ or $j^2$-skew-symmetric algebras
with $N$ generators the dimensions of their subspaces of given polynomial order are given by the following generating function:
\begin{equation}
H(t) = 1 + N t + N^2 t^2 + \frac{N (N-1)(N+1)}{3}\,t^3,
\end{equation}
where we include pure numbers (dimension $1$), the $N$ generators $\theta^A$ (or ${\bar{\theta}}^{\dot{B}}$), the $N^2$ independent quadratic combinations
$\theta^A \theta^B$ and $N(N-1)(N+1)/3$ products of three generators $\theta^A \theta^B \theta^C$.
\end{itemize}
%It is easy to see that all higher-order monomials starting from $4$-th power must identically vanish if associativity holds:
%\begin{equation}
%\theta^A \theta^B \theta^C \theta^D = j \, \theta^B \theta^C \theta^A \theta^D = j^2 \theta^B \theta^A \theta^D \theta^C =
%j^3 \theta^A \theta^D \theta^B \theta^C = j^4 \theta^A \theta^B \theta^C \theta^D.
%\end{equation}
%As $j^4 = j \neq 1$, the expression $\theta^a \theta^B \theta^C \theta^D$ must identically vanish.
\noindent
The above four cubic generalization of Grassmann algebra are represented in the following diagram, in which all the arrows are surjective homomorphisms.
\vskip 0.2cm
\centerline{$\Lambda_0$}
\centerline{$\downarrow$}
\centerline{$\Lambda_1$}
\centerline{$ \swarrow \hskip 0.4cm \searrow$}
\centerline{$ \Lambda \hskip 1.5cm {\bar{\Lambda}}$}
\vskip 0.2cm
\section{Examples of $Z_3$-graded ternary algebras}

\subsection{The $Z_3$-graded analogue of Grassman algebra}

Let us introduce $N$ generators spanning a linear space over complex numbers,
satisfying the following cubic relations \cite{Kerner1991}, \cite{Kerner1992}:

\begin{equation}
\theta^A \theta^B \theta^C = j \, \theta^B \theta^C \theta^A = j^2 \, \theta^C \theta^A \theta^B,
\label{ternary1}
\end{equation}
with $j = e^{2 i \pi/3}$, the primitive root of $1$. We have $1+j+j^2 = 0$ \; \; and ${\bar{j}} = j^2$.

Let us denote the algebra spanned by the  $\theta^A$ generators by ${\bf{\cal{A}}}$ \cite{Kerner1991}, \cite{Kerner1992}.

We shall also introduce a similar set of {\it conjugate} generators,
 ${\bar{\theta}}^{\dot{A}}$,
$\dot{A}, \dot{B},... = 1,2,...,N$, satisfying similar condition with $j^2$ replacing $j$:

\begin{equation}
{\bar{\theta}}^{\dot{A}} {\bar{\theta}}^{\dot{B}} {\bar{\theta}}^{\dot{C}} =
j^2 \, {\bar{\theta}}^{\dot{B}} {\bar{\theta}}^{\dot{C}} {\bar{\theta}}^{\dot{A}}
= j \, {\bar{\theta}}^{\dot{C}} {\bar{\theta}}^{\dot{A}} {\bar{\theta}}^{\dot{B}},
\label{ternary2}
\end{equation}

Let us denote this algebra by  ${\bar{\cal{A}}}$.

We shall endow the algebra ${\cal{A}} \oplus {\bar{\cal{A}}}$ with a natural  $Z_3$ grading, considering the generators  $\theta^A$
as grade $1$ elements, their conjugates ${\bar{\theta}}^{\dot{A}}$ being of grade $2$.

The grades add up modulo $3$, so that the products
 $\theta^{A} \theta^{B}$ span a linear
subspace of grade $2$, and the cubic products
$ \theta^A \theta^B \theta^C$ being of grade $0$.
Similarly, all quadratic expressions in conjugate generators,
${\bar{\theta}}^{\dot{A}} {\bar{\theta}}^{\dot{B}}$
are of  grade $2 + 2 = 4_{mod \, 3} = 1$, whereas their cubic products are again of grade $0$,
like the cubic products of $\theta^A$'s. \cite{Kerner2001}

Combined with the associativity, these cubic relations impose finite dimension on the
algebra generated by the $Z_3$ graded generators. As a matter of fact, cubic expressions are the
highest order that does not vanish identically. The proof is immediate:
\begin{equation}
\theta^A \theta^B \theta^C \theta^D = j \, \theta^B \theta^C \theta^A \theta^D =
j^2 \, \theta^B \theta^A \theta^D \theta^C
= j^3 \, \theta^A \theta^D \theta^B \theta^C =
j^4 \, \theta^A \theta^B \theta^C \theta^D,
\label{quartic1}
\end{equation}
and because
$j^4 = j \neq 1$, the only solution is $\theta^A \theta^B \theta^C \theta^D = 0.$
%\label{quartic2}
%\end{equation}

\subsection{The $Z_3$ graded differential forms}

Instead of the usual exterior differential operator satisfying
$d^2 = 0 $ , we can postulate its $Z_3$-graded generalization satisfying
$$d^2 \neq 0, \; \; \;   d^3 f = 0 $$
The first differential of a smooth function $f (x^i)$ is as usual
$$df = \partial_i f \, dx^i, $$
 whereas the second differential is formally
$$d^2 f = ( \partial_k \partial_i f ) \; dx^k dx^i + (\partial_i f ) \; d^2 x^i $$

 We shall attribute the grade $1$ to the $1$-forms
$d x^i, \; (i,j,k = 1,2,...N)$ , and  grade $2$  to the forms
$d^2 x^i, \; (i,j,k = 1,2,...N) $ ; under associative multiplication of these forms the grades
add up {\it modulo} $3$
$${\rm grade} ( \omega   \; \theta) = {\rm grade} (\omega) + {\rm grade} (\theta) \; (modulo 3).$$
 The $Z_3$-graded differential operator $d$
has the following property, compatible with grading we have chosen:
$$d (\omega \; \theta) = ( d \omega ) \, \theta + j^{{\rm grade}_{\omega}} \, \omega \, d \theta. $$

$$d^2 f = (\partial_i \partial_k f) d x^i d x^k + (\partial_i f ) \, d^2 x^i,$$

$$d^3 f = (\partial_m \partial_i \partial_k f ) d x^m d x^i d x^k +
( \partial_i \partial_k f ) d^2 x^i d x^k $$
$$+ j \; ( \partial_i \partial_k f ) d x^i d^2 x^k +
(\partial_k \partial_i f) d x^k d^2 x^i + (\partial_i f ) \, d^3 x^i .$$
equivalent with
$$d^3 f = (\partial_m \partial_i \partial_k f ) d x^m d x^i d x^k + ( \partial_i \partial_k f) [ d^2 x^k d x^i - j^2 \, d x^i d^2 x^k ]
+ (\partial_i f ) \, d^3 x^i.$$
 Consequently, assuming that $d^3 x^k = 0$ and $d^3 f = 0$, to make
the remaining terms vanish we must impose the following commutation relations
on the products of forms:
$$d x^i d x^k d x^m = j \, d x^k d x^m d x^i, \; \; \; \; \;
d x^i d^2 x^k = j \, d^2 x^k d x^i, $$
therefore
$$d^2 x^k d x^i = j^2 \, d x^i d^2 x^k$$

As in the case of the abstract $Z_3$-graded Grassmann algebra, the fourth order
expressions must vanish due to the associativity of the product:
$$dx^i dx^k dx^l dx^m =0.$$
Consequently, we shall assume that also
$$d^2 x^i d^2 x^k = 0.$$
This completes the construction of algebra of $Z_3$-graded exterior forms.

\subsection{Ternary Clifford algebra}
\indent
Let us introduce the following three $3 \times 3$ matrices:
\begin{equation}
Q_1 = \begin{pmatrix} 0 & 1 & 0 \cr 0 & 0 & j \cr j^2 & 0 & 0 \end{pmatrix}, \; \; \;
Q_2 = \begin{pmatrix} 0 & j & 0 \cr 0 & 0 & 1 \cr j^2 & 0 & 0 \end{pmatrix}, \; \; \;
Q_3 = \begin{pmatrix} 0 & 1 & 0 \cr 0 & 0 & 1 \cr 1 & 0 & 0 \end{pmatrix},
\label{threeQ}
\end{equation}
and their hermitian conjugates
\begin{equation}
Q^{\dagger}_1 = \begin{pmatrix}  0 & 0 & j \cr 1 & 0 & 0 \cr 0 & j^2 & 0 \end{pmatrix}, \; \; \;
Q^{\dagger}_2 = \begin{pmatrix} 0 & 0 & j \cr j^2 & 0 & 0 \cr 0 & 1 & 0 \end{pmatrix}, \; \; \;
Q^{\dagger}_3 = \begin{pmatrix} 0 & 0 & 1 \cr 1 & 0 & 0 \cr 0 & 1 & 0 \end{pmatrix}.
\label{threeQbar}
\end{equation}
These matrices can be allowed natural $Z_3$ grading,
\begin{equation}
{\rm grade} (Q_k) = 1, \; \; \; {\rm grade} (Q^{\dagger}_k) = 2,
\label{gradeQ}
\end{equation}
The above matrices span a very interesting ternary algebra. Out of three independent $Z_3$-graded ternary
combinations, only one leads to a non-vanishing result. One can check without much effort that both $j$ and $j^2$ skew
ternary commutators do vanish:
$$
\{ Q_1, Q_2, Q_3 \}_j = Q_1 Q_2 Q_3 + j Q_2 Q_3 Q_1 + j^2 Q_3 Q_1 Q_2 = 0, $$
$$\{ Q_1, Q_2, Q_3 \}_{j^2} = Q_1 Q_2 Q_3 + j^2 Q_2 Q_3 Q_1 + j Q_3 Q_1 Q_2 = 0,$$
and similarly for the odd permutation, $Q_2 Q_1 Q_3$.
On the contrary, the totally symmetric combination does not vanish; it is proportional to the $3 \times 3$ identity matrix ${\bf 1}$:
\begin{equation}
 Q_a Q_b Q_c + Q_b Q_c Q_a + Q_c Q_a Q_b = 3\,\eta_{abc} \, {\bf 1}, \; \; \; a,b,... = 1,2,3.
\label{anticom}
\end{equation}
with $\eta_{abc}$ given by the following non-zero components:
\begin{equation}
\eta_{111} = \eta_{222} = \eta_{333} = 1, \; \; \; \eta_{123} = \eta_{231} = \eta_{312} = 1, \; \; \;
\eta_{213} = \eta_{321} = \eta_{132} = j^2.
\label{defeta}
\end{equation}
all other components vanishing. The relation ({\ref{anticom}) may serve as the definition of {\it ternary Clifford algebra}.

Another set of three matrices is formed by the hermitian conjugates of $Q_a$, which we shall endow with dotted indeces ${\dot{a}}, {\dot{b}},...=1,2,3$:
%\begin{equation}
$Q_{{\dot{a}}} = Q_a^{\dagger}$
%\label{defqbar}
%\end{equation}
satisfying conjugate identities
 \begin{equation}
 Q_{\dot{a}} Q_{\dot{b}} Q_{\dot{c}} + Q_{\dot{b}} Q_{\dot{c}} Q_{\dot{a}} + Q_{\dot{c}} Q_{\dot{a}} Q_{\dot{b}}
 = 3\,\eta_{{\dot{a}}{\dot{b}}{\dot{c}}} \, {\bf 1}, \; \; \; {\dot{a}}, {\dot{b}},... = 1,2,3.
\label{anticomdot}
\end{equation}
with $\eta_{{\dot{a}}{\dot{b}}{\dot{c}}} = {\bar{\eta}}_{abc}$.

It is obvious that any similarity transformation of the generators $Q_a$ will keep the ternary anti-commutator (\ref{anticom})
invariant. As a matter of fact, if we define ${\tilde{Q}}_b = P^{-1} Q_b P$, with $P$ a non-singular $3 \times 3$ matrix,
 the new set of generators will satisfy the same ternary relations, because
$${\tilde{Q}}_a {\tilde{Q}}_b {\tilde{Q}}_c = P^{-1} Q_a P P^{-1} Q_b P P^{-1} Q_c P = P^{-1} (Q_a Q_b Q_c) P,$$
and on the right-hand side we have the unit matrix which commutes with all other matrices, so that $P^{-1} \; {\bf 1} \; P  = {\bf 1}$.

\section{Generalized $Z_2 \times Z_3$-graded ternary algebra}

Let us suppose that we have binary and ternary skew-symmetric products defined by corresponding structure constants:
\begin{equation}
\xi^{\alpha} \xi^{\beta} = - \xi^{\beta} \xi^{\alpha}
\label{skewxi}
\end{equation}
\begin{equation}
\theta^A \theta^B \theta^C = j \; \theta^B \theta^C \theta^A
\label{thetaskew}
\end{equation}
The unifying ternary relation is of the type $\Lambda_0$, i.e.
\begin{equation}
X^i X^j X^k + X^j X^k X^i + X^k X^i X^j + X^k X^j X^i + X^j X^i X^k + X^i X^k X^j = 0.
\label{sixsum}
\end{equation}
It is obviously satisfied by both types of variables; the $\theta^A$'s by definition of the product,
for which at this stage the associativity property can be not decided yet; the product of
grassmannian  $\xi^{\alpha}$ variables (\ref{skewxi}) on the contrary, should be associative in order
to make the formula (\ref{sixsum}) applicable.

It can be added that the cubic constitutive relation (\ref{thetaskew}) satisfies a simpler condition
with cyclic permutations only,
$$\theta^A \theta^B \theta^C + \theta^B \theta^C \theta^A + \theta^C \theta^A \theta^B = 0,$$
but the cubic products of grassmannian variables are invariant under even (cyclic) permutations,
so that only the combination of all six permutations of $\xi^{\alpha} \xi^{\beta} \xi^{\gamma}$,
like in (\ref{sixsum}) will vanish.

Now, if we want to merge the two algebras into a common one, we must impose the general condition (\ref{sixsum})
to the mixed cubic products. These are of two types: $\theta^A \xi^{\alpha} \theta^{B}$ and $\xi^{\alpha} \theta^B \xi^{\beta}$,
with two $\theta$'s and one $\xi$, or with two $\xi$'s and one $\theta$. These identities, all like (\ref{sixsum}) should follow
from {\it binary} constitutive relations imposed on the {\it associative} products between one $\theta$ and one $\xi$ variable.

Let us suppose that one has
\begin{equation}
\xi^{\alpha} \theta^B = \omega \; \theta^B \xi^{\alpha} \; \; \; {\rm and \; consequently} \; \; \;
 \theta^A \xi^{\beta} = \omega^{-1} \; \xi^{\beta} \theta^A.
\label{mixedxitheta}
\end{equation}
A simple exercise leads to the conclusion that in order to satisfy the general condition (\ref{sixsum}), the
unknown factor $\omega$ must verify the equation $\omega + \omega^{-1} + 1=0$, or equivalently, $\omega + \omega^2 + \omega^3 = 0$. Indeed, we have, assuming
 the associativity:
$$\theta^A \xi^{\alpha} \theta^B = \omega^{-1} \; \xi^{\alpha} \theta^A \theta^B = \omega \; \theta^A \theta^B \xi^{\alpha},$$
$$\theta^B \xi^{\alpha} \theta^A = \omega^{-1} \; \xi^{\alpha} \theta^B \theta^A = \omega \; \theta^B \theta^A \xi^{\alpha}.$$
>From this we get, by transforming all the six products so that $\xi^{\alpha}$ should appear always in front of
the monomials:
$$\theta^A \xi^{\alpha} \theta^B = \omega^{-1} \; \xi^{\alpha} \theta^A \theta^B, \; \; \theta^A \theta^B \xi^{\alpha} = \omega^{-2} \; \xi^{\alpha} \theta^A \theta^B,$$
$$\theta^B \xi^{\alpha} \theta^A = \omega^{-1} \; \xi^{\alpha} \theta^B \theta^A, \; \;  \theta^B \theta^A \xi^{\alpha} =    \omega^{-2} \;  \xi^{\alpha} \theta^B \theta^A.$$
Adding up all permutatuons, even (cyclic) and odd alike, we get the following result:
$$ \theta^A \xi^{\alpha} \theta^B + \xi^{\alpha} \theta^B \theta^A + \theta^B \theta^A \xi^{\alpha} +
\theta^B \xi^{\alpha} \theta^A + \xi^{\alpha} \theta^A \theta^B + \theta^A \theta^B \xi^{\alpha} =$$
\begin{equation}
(1 + \omega + \omega^{-1}) \; \xi^{\alpha} \theta^A \theta^B + (1 + \omega + \omega^{-1}) \; \xi^{\alpha} \theta^B \theta^A.
\label{twocycles}
\end{equation}
The expression in (\ref{twocycles}) will identically vanish if $\omega = j = e^{\frac{2 \pi i}{3}}$ (or $j^2$, which satisfies the same relation
$j+j^2 +1 =0.$

The second type of cubic monomials, $\xi^{\alpha} \theta^B \xi^{\beta}$, satisfies the identity
\begin{equation}
\xi^{\alpha} \theta^B \xi^{\delta} + \theta^B \xi^{\delta} \xi^{\alpha} + \xi^{\delta} \xi^{\alpha} \theta^B +
\xi^{\delta} \theta^B \xi^{\alpha} + \theta^B \xi^{\alpha} \xi^{\delta} + \xi^{\alpha} \xi^{\delta} \theta^B = 0
\label{xithetaxi}
\end{equation}
no matter what the value of $\omega$ is chosen in the constitutive relation (\ref{mixedxitheta}), the antisymmetry of the product of two
$\xi$'s suffices. As a matter of fact, because we have $\xi^{\alpha} \xi^{\delta} = - \xi^{\delta} \xi^{\alpha}$, in the formula
(\ref{xithetaxi}) the second term cancels the fifth term, and the third termis cancelled by the sixth one. What remains is the sum of the
first and the fourth terms:
$$\xi^{\alpha} \theta^B \xi^{\delta} + \xi^{\delta} \theta^B \xi^{\alpha}.$$
Now we can transform both terms so as to put the factor $\theta$ in front; this will give
\begin{equation}
\xi^{\alpha} \theta^B \xi^{\delta} + \xi^{\delta} \theta^B \xi^{\alpha} =
\omega \theta^B \xi^{\alpha} \xi^{\delta} + \omega \theta^B  \xi^{\delta} \xi^{\alpha} = 0
\end{equation}
because of the anti-symmetry of the product between the two $\xi$'s.

This completes the construction of the $Z_3 \times Z_2$-graded extension of Grassman algebra.

The existence of {\it two } cubic roots of unity, $j$ and $j^2$, suggests that one can extend the above algebraic construction
by introducing a set of {\it conjugate} generators, denoted for convenience with a bar and with dotted indeces, satisfying
conjugate ternary constitutive relation (\ref{ternary2}). The unifying condition of vanishing of the sum of all
permutations (algebra of $\Lambda_0$-type) will be automatically satisfied.

But now we have to extend this condition to the triple products of the type $\theta^A {\bar{\theta}}^{\dot{B}} \theta^C$ and
${\bar{\theta}}^{\dot{A}} \theta^B {\bar{\theta}}^{\dot{C}}$. This will be achieved if we impose the obvious
condition, similar to the one proposed already fo binary combinations $\xi \theta$:
\begin{equation}
\theta^A {\bar{\theta}}^{\dot{B}} = j {\bar{\theta}}^{\dot{B}} \theta^A, \; \; \;
{\bar{\theta}}^{\dot{B}} \theta^A = j^2 \; \theta^A {\bar{\theta}}^{\dot{B}}
\label{thetathetabar}
\end{equation}
The proof of the validity of the condition (\ref{sixsum}) for the above combinations is exactly the same as for the
triple products $\xi^{\alpha} \theta^B \xi^{\gamma}$ and $\theta^A \xi^{\delta} \theta^B$.

We have also to impose commutation relations on the mixed products of the type
$$\xi^{\alpha} {\bar{\theta}}^{\dot{B}} \xi^{\beta} \; \; \; {\rm and} \; \; \; {\bar{\theta}}^{\dot{B}} \xi^{\beta} {\bar{\theta}}^{\dot{C}}.$$
It is easy to see that like in the former case, it is enough to impose the commutation rule similar to the former one with $\theta$'s,
namely
\begin{equation}
{\xi}^{\alpha} {\bar{\theta}}^{\dot{B}} = j^2 \; {\bar{\theta}}^{\dot{B}} {\xi}^{\alpha}
\label{xithetabar}
\end{equation}
Although we could stop at this point the extension of our algebra, for the sake of symmetry it seems useful to
introduce the new set of conjugate variables ${\bar{\xi}}^{\dot{\alpha}}$ of the $Z_2$-graded type. We shall suppose
that they anti-commute, like the $\xi^{\beta}$'s, and not only between themselves, but also with their conjugates,
which means that we assume
\begin{equation}
{\bar{\xi}}^{\dot{\alpha}} {\bar{\xi}}^{\dot{\beta}} = - {\bar{\xi}}^{\dot{\beta}} {\bar{\xi}}^{\dot{\alpha}},\;\;\;\;
 {{\xi}^{\alpha}} {\bar{\xi}}^{\dot{\beta}} = - {\bar{\xi}}^{\dot{\beta}} {{\xi}^{\alpha}}.
\label{antoicomxibar}
\end{equation}
This ensures that the condition (\ref{sixsum}) will be satisfied by any ternary combination of the $Z_2$-graded generators,
including the mixed ones like
$${\bar{\xi}}^{\dot{\alpha}} \xi^{\beta} {\bar{\xi}}^{\dot{\delta}} \; \; \; {\rm or} \; \; \;  \xi^{\beta} {\bar{\xi}}^{\dot{\alpha}} \xi^{\gamma}.$$
The dimensions of classical Grassmann algebras with $n$ generators are well known: they are equal to $2^n$, with subspaces spanned by the
products of $k$ generators having the dimension $C^n_k = n!/(n-k)!k!$. With $2n$ anticommuting generators, $\xi^{\alpha}$ and ${\bar{\xi}}^{\dot{\beta}}$
we shall have the dimension of the corresponding Grassmann algebra equal to $2^{2n}$.

It is also quite easy to determine the dimension of the $Z_3$-graded generalizations of Grassmann algebras constructed above
(see, e.g. in \cite{Kerner1992}, \cite{Abramov1997}, \cite{Kerner1997} ). The $Z_3$-graded algebra with $N$ generators $\theta^A$
has the total dimension $N + N^2 + (N^3 - N)/3 = (N^3 + N^2 + 2N)/3$. The conjugate algebra, with the same number of generators,
has the identical dimension. However, the dimension of the extended algebra unifying both these algebras is not equal to
the square of the dimension of one of them because of the extra conditions on the mixed products between the generators and their
conjugates, $\theta^A {\bar{\theta}}^{\dot{B}} = {\bar{\theta}}^{\dot{B}} \theta^A.$

\section{Two distinct gradings: $Z_3 \times Z_2$ versus $Z_6$}

The $Z_2$-grading of ordinary (binary) algebras is well known and widely studied and applied (e.g. in the super-symmetric
field theories in Physics).  The Grassmann algebra is perhaps the oldest and the best known example of a $Z_2$-graded structure.
Other gradings are much less popular. The $Z_3$-grading was introduced and studied in the paper \cite{Kerner1991};
the $Z_N$ grading was discussed in \cite{MDVRK1996} and \cite{MDV1998}.

In the case of ternary algebras of type $\Lambda_1$ or $\Lambda_2$, the grade $1$ is attributed to the generators $\theta^A$
and the grade $2$ to the conjugate generators ${\bar{\theta}}^{\dot{B}}$. Consequently, their products acquire the grade
which is the sum of grades of the factors modulo $3$. When we consider an algebra including a ternary $Z_3$-graded subalgebra
 and a binary $Z_2$-graded one, we can quite naturally introduce a combination of the two gradings considered as a pair
of two numbers, say $(a, \lambda)$, with $a = 0,1,2$ representing the $Z_3$-grade, and $\lambda = 0, 1$ representing the
$Z_2$ grade, $\lambda = 0,1$. The first grades add up modulo $3$, the second grades add up modulo $2$.
The six possible combined grades are then
\begin{equation}
(0,0), \; \; \; (1, 0), \; \; \; (2, 0), \; \; \; (0, 1), \; \; \; (1,1) \; \; \; {\rm and} \; \; \; (2,1).
\label{combgrade}
\end{equation}
To add up two of the combined grades amounts to adding up their first entries modulo $3$, and their second entries
modulo $2$. Thus, we have
$$(2,1) + (1,1) = (3,2) \simeq (0,0), \; \; \; {\rm or} \; \; \; ((2,1) + (1, 0) = (3,1) \simeq (0,1), \; \; \;
{\rm and \; so \; forth.}$$
It is well known that the cartesian product of two cyclic groups $Z_N \times Z_n$, $N$ and $n$ being two prime numbers,
 is the cyclic group $Z_{Nn}$ corresponding to the product of those prime numbers. This means that there is an isomorphism
between the cyclic group $Z_6$, generated by the {\it sixth} primitive root of unity, $q^6 = 1$, satisfying the equation
$$q + q^2 + q^3 + q^4 + q^5 + q^6 = 0.$$
This group can be represented on the complex plane, with $q =e^{\frac{2 \pi i}{6}}$, as shown on the diagram below:

%%%%%%%%%%%%%%%%%%%%%%%%%%%%%% Figure begins
%\begin{comment}
\begin{figure}[hbt]
\centering
\includegraphics[width=5cm, height=5cm]{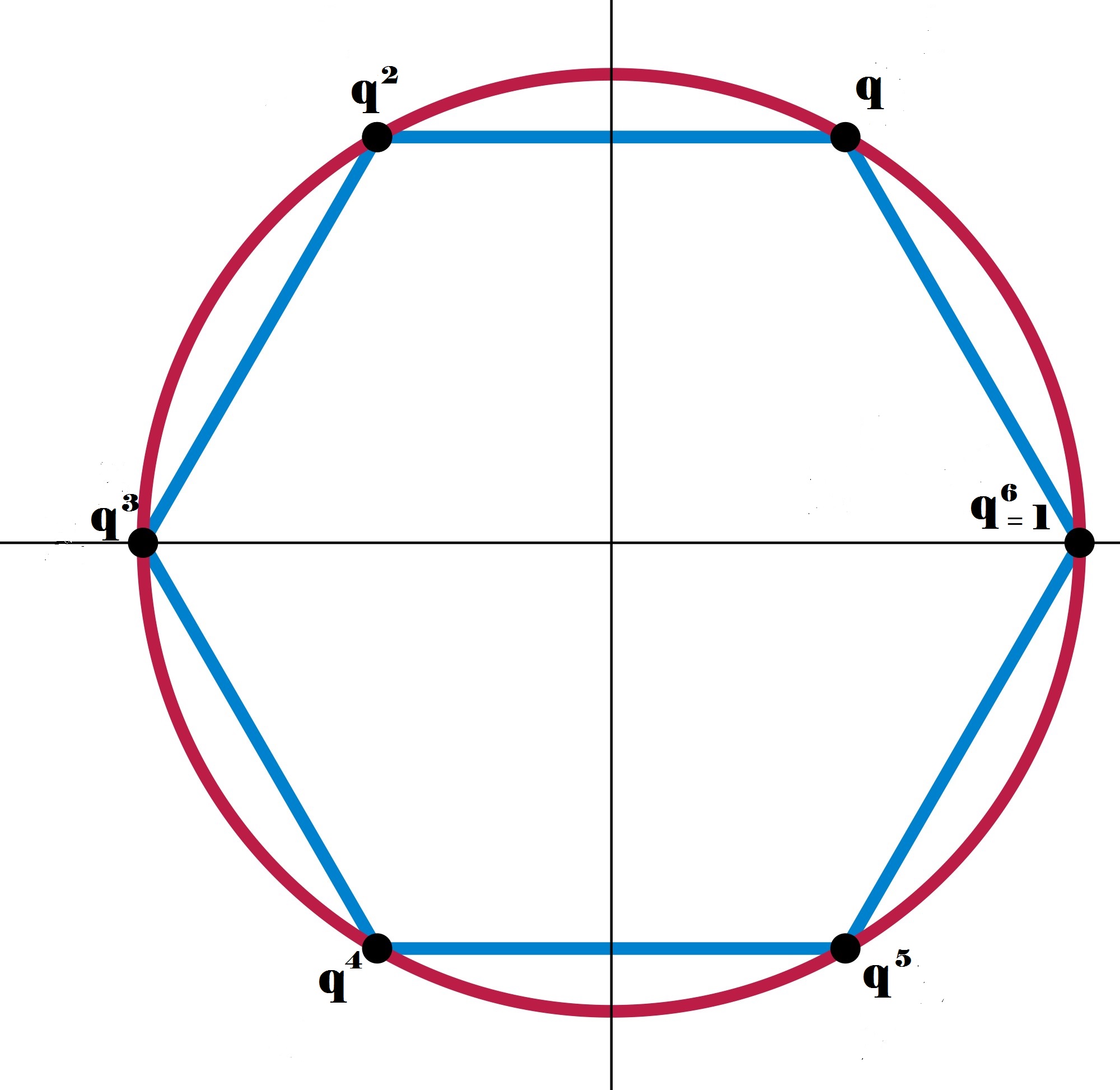}
\hskip 0.9cm
\includegraphics[width=5cm, height=5cm]{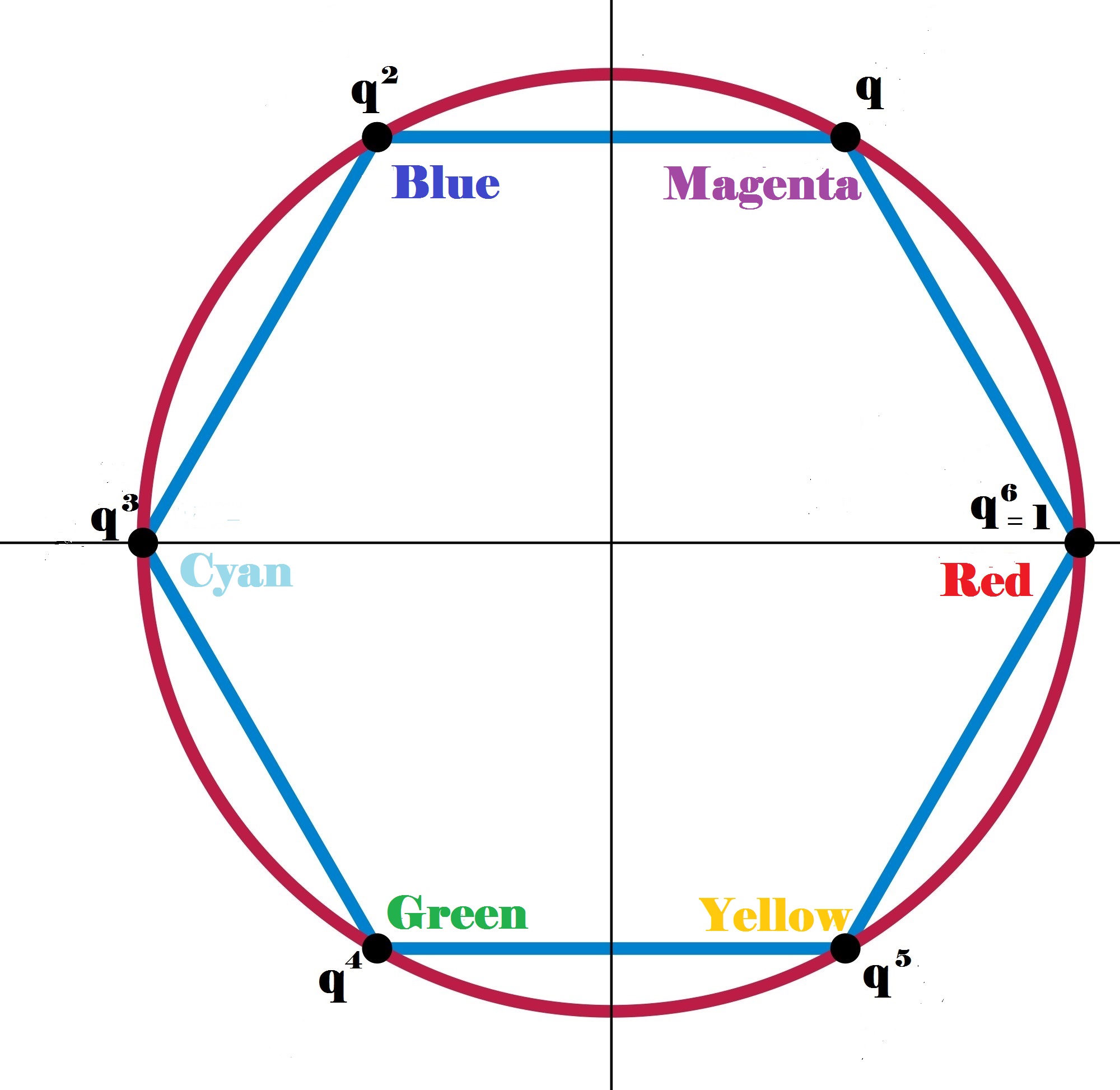}
\caption{{\small Representation of the cyclic group $Z_6$ in the complex plane with three colors and three "anti-colors"
 attributed to even and odd powers of $q$, accordingly with colors attributed in Quantum Chromodynamics to quatks and to anti-quarks.}}
\label{fig:Cyclic_q6}
\end{figure}
%\end{comment}
%%%%%%%%%%%%%%%%%%%%%%%%%%%%% Figure ends

%%%%%%%%%%%%%%%%%%%%%%%%%%%%%% Figure begins
%\begin{comment}
%\begin{figure}[hbt]
%\centering
%\includegraphics[width=6cm, height=6cm]{Z6_system.jpg}
%\caption{{\small Representation of the cyclic group $Z_6$ in the complex plane.}}
%\label{fig:sixroots}
%\end{figure}
%\end{comment}
%%%%%%%%%%%%%%%%%%%%%%%%%%%%% Figure ends
The elements of the group $Z_6$ represented by complex numbers  multiply modulo $6$, e. g. $q^4 \cdot q^5 = q^9 \simeq q^3,$ etc.
The six elements of $Z_6$ can be put in the one-to-one correspondence with the pairs defining six elements of $Z_3 \times Z_2$
according to the following scheme:
\begin{equation}
(0,0) \simeq q^0 = 1, \; \; (2,1) \simeq q, \; \; (1,0) \simeq q^2, \; \; (0,1) \simeq q^3, \; \; (2,0) \simeq q^4, \; \; (1,1) \simeq q^5.
\label{Z6Z3Z2}
\end{equation}
The same result can be obtained directly using the representations of $Z_3$ and $Z_2$ in the complex plane. Taken separately, each
of these cyclic groups is generated by one non-trivial element, the third root of unity $j = e^{\frac{2 \pi i}{3}}$ for $Z_3$
and $-1 = e^{\pi i}$ for $Z_2$. It is enough to multiply these complex numbers and take their different powers in order to
get all the six elements of the cyclic group $Z_6$. One easily identifies then
$$-j^2 = q, \; \; j= q^2, \; \;  - 1 = q^3, \; \; j^2 = q^4, \; \; -j = q^5, \; \; \; 1 = q^6.$$
The colors attributed to the powers of the complex generator $q$ can be used to modelize the exclusion principle used in 
Quantum Chromodynamics, where exclusively the "white" combinations of three quarks and three anti-quarks, as well as
 the "white" quark-anti-quark pairs are declared observable. Replacing the word "white" by $0$, we see that there
are {\it two} vanishing linear combinations of {\it three} powers of $q$, and {\it three} pairs of powers of $q$ that 
are also equal to zero. Indeed, we have:
\begin{equation}
q^2 + q^4 + q^6 = j + j^2 + 1 =0, \; \; {\rm and} \; \; q + q^3 + q^5 = -j^2 -1 -j = 0,
\label{twotriq}
\end{equation}
as well as 
\begin{equation}
q + q^4 = 0, \; \; \; q^2 + q^5 = 0, \; \; q^3 + q^6 = 0.
\label{threetwoq}
\end{equation}
 The $Z_6$-grading should unite both $Z_2$ and $Z_3$ gradings, reproducing their essential properties.
Obviously, the $Z_3$ subgroup is formed by the elements $1, \; q^2$ and $q^4$, while the $Z_2$ subgroup
is formed by the elements $1$ and $q^3 = -1$.
In what follows, we shall see
that the associativity imposes many restrictions which can be postponed in the case of non-associative ternary structures.

The natural choice for the $Z_3$-graded algebra with cubic relations was to attribute the grade $1$ to the generators $\theta^A$,
and grade $2$ to their conjugates ${\bar{\theta}}^{\dot{B}}$. All other expressions formed by products and powers of those
 got the well defined grade, the sum of the grades of factors modulo $3$. In a simple Cartesian product of two algebras,
a $Z_3$-graded with a $Z_2$-graded one, the generators of the latter will be given grade $1$, and their products
will get automatically the grade which is the sum of the grades of factors modulo $2$, which means that the all products
and powers of generators $\xi^{\alpha}$ will acquire grade $1$ or $0$ according to the number and character of
factors involved. The mixed products of the type $\theta^A \xi^{\beta}, \; \xi^{\beta} \theta^B \theta^C$, etc.
can be given the double $Z_3 \times Z_2$ grade according to (\ref{combgrade}). According to the isomorphism
defined by (\ref{Z6Z3Z2}), this is equivalent to a $Z_6$-grading of the product algebra.

As long as the algebra is supposed to be {\it homogeneous} in the sense that all the constitutive relations
contain exclusively terms of {\it one and the same type}, like in the extension of Grassmann algebra discussed
above, the supposed associativity does not impose any particular restrictions. However, this is not the case
if we consider the possibility of {\it non-homogeneous} constitutive equations, including terms of different nature,
but with the same $Z_6$-grade. The grading defined by (\ref{Z6Z3Z2}) suggests a possibility of extending the
constitutive relations by comparing terms of the type $\theta^A \theta^B \theta^C$, whose $Z_6$-grade is $3$,
to  the generators $\xi^{\alpha}$ having the same $Z_6$-grade. This will lead to the following constitutive
relations:
\begin{equation}
\theta^A \theta^B \theta^C = \rho^{ABC}_{\; \; \; \; \; \; \; \; \alpha} \; \xi^{\alpha} \; \; \; \; {\rm and}  \; \; \; \;
{\bar{\theta}}^{\dot{A}} {\bar{\theta}}^{\dot{B}} {\bar{\theta}}^{\dot{C}} =
{\bar{\rho}}^{{\dot{A}} {\dot{B}} {\dot{C}}}_{\; \; \; \; \; \; \; \; \dot{\alpha}} \; {\bar{\xi}}^{\dot{\alpha}}
\label{nonhomthetaxi}
\end{equation}
with the coefficients (structure constants) $\rho^{ABC}_{\; \; \; \; \; \; \; \;  \alpha}$ and
${\bar{\rho}}^{{\dot{A}} {\dot{B}} {\dot{C}}}_{\; \; \; \; \; \; \; \dot{\alpha}}$ displaying obvious symmetry properties mimicking
the properties of ternary products of $\theta$-generators with respect to cyclic permutations:
\begin{equation}
\rho^{ABC}_{\; \; \; \; \; \; \; \;  \alpha} = j \; \rho^{BCA}_{\; \; \; \; \; \; \; \; \alpha}
= j^2 \; \rho^{CAB}_{\; \; \; \; \; \; \; \; \alpha} \; \; \; {\rm and} \; \; \;
{\bar{\rho}}^{{\dot{A}} {\dot{B}} {\dot{C}}}_{\; \; \; \; \; \; \; \; \dot{\alpha}} =
j^2 \; {\bar{\rho}}^{{\dot{B}} {\dot{C}} {\dot{A}}}_{\; \; \; \; \; \; \; \; \dot{\alpha}} =
j \; {\bar{\rho}}^{{\dot{C}} {\dot{A}} {\dot{B}}}_{\; \; \; \; \; \; \; \; \dot{\alpha}}.
\label{rhoprop}
\end{equation}
If all products are supposed to be associative, then we see immediately that the products between $\theta$ and $\xi$ generators,
as well as those between ${\bar{\theta}}$ and ${\bar{\xi}}$  generators must vanish identically, because of the vanishing of
quatric products $\theta \theta \theta \theta =0$ and ${\bar{\theta}} {\bar{\theta}} {\bar{\theta}} {\bar{\theta}} = 0$.
This means that we must set
\begin{equation}
\theta^A \; \xi^{\beta} = 0, \; \; \; \xi^{\beta} \theta^A = 0, \; \; \; {\rm as \; well \; as} \; \; \;
{\bar{\theta}}^{\dot{B}} {\bar{\xi}}^{\dot{\alpha}} = 0, \; \; \; {\bar{\xi}}^{\dot{\alpha}} {\bar{\theta}}^{\dot{B}} = 0.
\label{thetaxicom}
\end{equation}
But now we want to unite the two gradings into a unique common one. Let us start by defining a ternary product of generators,
not necessarily derived from an ordinary associative algebra. We shall just suppose the existence of ternary product of
generators, displaying the $j$-skew symmetry property:
\begin{equation}
\{ \theta^A,  \theta^B , \theta^C \} = j  \{ \theta^B, \theta^C, \theta^A \} = j^2 \{ \theta^C,  \theta^A, \theta^B \}.
\label{terprod1}
\end{equation}
and similarly, for the conjugate generators,
\begin{equation}
\{ {\bar{\theta}}^{\dot{A}}, {\bar{\theta}}^{\dot{B}}, {\bar{\theta}}^{\dot{C}} \} =
 j^2 \; \{ {\bar{\theta}}^{\dot{B}}, {\bar{\theta}}^{\dot{C}}, {\bar{\theta}}^{\dot{A}} \} =
j \; \{ {\bar{\theta}}^{\dot{C}}, {\bar{\theta}}^{\dot{A}}, {\bar{\theta}}^{\dot{B}} \}.
\label{terprod2}
\end{equation}
Let us attribute the $Z_6$-grade $1$ to the generators $\theta^A$. Then it is logical to attribute the $Z_6$ grade $5$
to the conjugate generators ${\bar{\theta}}^{\dot{B}}$, so that mixed  products $\theta^A {\bar{\theta}}^{\dot{B}}$ would
be of $Z_6$ grade $0$. Ternary products (\ref{terprod1}) are of grade $3$, and ternary products of conjugate generators
(\ref{terprod2}) are also of grade $3$, because $5+5+5=15$, and $15 \; modulo \;  6 = 3$. But we have also $q^3 = -1$, which
is the generator of the $Z_2$-subalgebra of $Z_6$. Therefore we should attribute the $Z_6$-grade $3$ to both kinds of
the anti-commuting variables, $\xi^{\alpha}$ and ${\bar{\xi}}^{\dot{\beta}}$, because we can write their constitutive relations
using the root $q$ as follows:
\begin{equation}
\xi^{\alpha} \xi^{\beta} = - \xi^{\beta} \xi^{\alpha}= q^3 \; \xi^{\beta} \xi^{\alpha}, \; \; \; \;
{\bar{\xi}}^{\dot{\alpha}} {\bar{\xi}}^{\dot{\beta}} = - {\bar{\xi}}^{\dot{\beta}} {\bar{\xi}}^{\dot{\alpha}} =
q^3 \;{\bar{\xi}}^{\dot{\beta}} {\bar{\xi}}^{\dot{\alpha}}, \; \; \; \;
\xi^{\alpha} {\bar{\xi}}^{\dot{\beta}} = - {\bar{\xi}}^{\dot{\beta}} \xi^{\alpha} = q^3 \; {\bar{\xi}}^{\dot{\beta}} \xi^{\alpha} ,
\label{xixibarcom}
\end{equation}
On the other hand, the expressions containing products of $\theta$ with ${\bar{\xi}}$ and ${\bar{\theta}} $ with $\xi$:
$$\theta^A {\bar{\xi}}^{\dot{\alpha}} \; \; \; {\rm and} \; \; \; {\bar{\theta}}^{\dot{B}} \xi^{\beta}$$
The first expression has the $Z_6$-grade $1+3 =4$, and the second product has the $Z_6$-grade $5+3 = 8$ {\it modulo} $6$ $=2$.
Other products endowed with the same grade in our associative $Z_6$-grade algebra are ${\bar{\theta}}^{\dot{A}} \, {\bar{\theta}}^{\dot{B}}$
(grade $4$, because $5 + 5 = 10$ \; {\it modulo} \, $6$ $=4$, and $\theta^{A} \theta^{B}$ of grade $2$, because $1+1=2$).

This suggests that the following non-homogeneous constitutive relations can be proposed:
\begin{equation}
\theta^A {\bar{\xi}}^{\dot{\alpha}} = f^{A {\dot{\alpha}}}_{\; \; \; \; \; {\dot{C}} {\dot{D}}} \; {\bar{\theta}}^{\dot{C}}  {\bar{\theta}}^{\dot{D}},
\; \; \; {\rm and} \; \; \;
{\bar{\theta}}^{\dot{A}} \xi^{\alpha} = {\bar{f}}^{{\dot{A}} \alpha}_{\; \; \; \; \; C D} \; \theta^C  \theta^D,
\label{relthetaxi}
\end{equation}
where the coefficients should display the symmetry properties contravariant to those of the generators themselves, which means that we should have
\begin{equation}
f^{A {\dot{\alpha}}}_{\; \; \; \; \; {\dot{C}} {\dot{D}}} = j^2 \; f^{{\dot{\alpha}} A}_{\; \; \; \; \; {\dot{C}} {\dot{D}}} \; \; \;
{\rm and} \; \; \;
{\bar{f}}^{{\dot{A}} \alpha}_{\; \; \; \; \; C D} = j \; {\bar{f}}^{ \alpha {\dot{A}} }_{\; \; \; \; \; C D}
\end{equation}

\section{Invariance groups of $Z_2$ and $Z_3$ skew-symmetric algebras}

Before discussing the merger of a $Z_2$-graded algebra with a $Z_3$-graded algebra, let us explore the invariance properties of
each type separately. First, let us ask what kind of linear transformations preserves the ordinary Grassmann algebra spanned
by anti-commuting generators $\xi^{\alpha} \; \; (\alpha = 1, 2,...n).$ Any linearly independent $n$ combinations of anti-commuting variables
$\xi^{\alpha}$ will span another anti-commuting basis: indeed, if $\eta^{\alpha'} = S^{\alpha'}_{\beta} \xi^{\beta}$, and take
on purely numerical values, i.e. do commute with all other generators, then we can write
\begin{equation}
 \eta^{\alpha'} \eta^{\delta'} = S^{\alpha'}_{\alpha} \xi^{\alpha} \; S^{\delta'}_{\delta} \xi^{\delta} =
S^{\alpha'}_{\alpha}  \; S^{\delta'}_{\delta} \; \xi^{\alpha} \xi^{\delta} = - S^{\alpha'}_{\alpha}  \; S^{\delta'}_{\delta} \; \xi^{\delta} \xi^{\alpha} 
= -  S^{\delta'}_{\delta} \; \xi^{\delta}  S^{\alpha'}_{\alpha} \;\xi^{\alpha} = - \eta^{\delta'} \eta^{\alpha'}
\label{etaxi}
\end{equation}
Let us consider the simplest case of a $Z_2$-graded algebra spanned by two generators $\xi^{\alpha}, \; \; \alpha, \beta =1,2.$
The anti-commutation property can be encoded in the invariant $2$-form $\varepsilon_{\alpha \beta} $.
%= - \varepsilon{\beta \alpha}$.
We can obviously write
$$\varepsilon_{\alpha \beta} \xi^{\alpha} \xi^{\beta} = \varepsilon_{\beta \alpha} \xi^{\beta} \xi^{\alpha} =
 - \varepsilon_{\beta \alpha} \xi^{\alpha} \xi^{\beta},$$
from which we conclude that  $\varepsilon_{\alpha \beta} = - \varepsilon_{\beta \alpha}$.
We can choose the basis in which
$$\varepsilon_{11} = 0, \; \; \; \varepsilon_{22} = 0, \; \; \; \varepsilon_{12} = - \varepsilon_{21} = 1.$$
After a change of basis, $\xi^{\beta} \rightarrow S^{\alpha'}_{\beta} \xi^{\beta} = \eta^{\alpha'}$ the $2$-form
$\varepsilon_{\alpha \beta}$, as any tensor,  also undergoes the inverse transformation:
$$ S^{\alpha}_{\alpha'} S^{\beta}_{\beta'} \, \varepsilon_{\alpha \beta},$$
with $S^{\alpha}_{\beta'}$ the inverse matrix of the matrix $S^{\beta'}_{\beta}$.
Whatever non-singular linear transformation $S^{\beta}_{\beta'}$ is chosen, the new components $\varepsilon_{\alpha' \beta'}$
remain anti-symmetric, but they have not necessarily the same values as those of $\varepsilon_{\alpha \beta}$
However, if we require that also in new basis
$$\varepsilon_{1'1'} = 0, \; \; \; \varepsilon_{2'2'} = 0, \; \; \; \varepsilon_{1'2'} = - \varepsilon_{2'1'} = 1,$$
then it is easy to show that this imposes extra condition on the $2 \times 2$ matrix $S^{\alpha'}_{\beta}$, namely that
det $S= 1$. This defines the $SL(2, {\bf C})$ group as the group of invariance of the subalgebra spanned by two anti-commuting
 generators $\xi^{\alpha}$, $\alpha, \beta, ... = 1,2$.

Now let us turn to the invariance properties of the ternary subalgebra spanned by two generators $\theta^1, \; \theta^2$,
satisfying homogeneous cubic $j$-anticummutation relations $\theta^A \theta^B \theta^C = q^2 \; \theta^B \theta^C \theta^A$,
and their conjugate counterparts ${\bar{\theta}}^{\dot{1}}, \; {\bar{\theta}}^{\dot{2}}$ satisfying homogeneous cubic
$j^2$ anti-commutation relations ${\bar{\theta}}^{\dot{A}} {\bar{\theta}}^{\dot{B}} {\bar{\theta}}^{\dot{C}} 
= j^2 \; {\bar{\theta}}^{\dot{B}} {\bar{\theta}}^{\dot{C}} {\bar{\theta}}^{\dot{A}}$

We shall also impose {\it binary} constitutive relations between the generators $\theta^A$ and their
conjugate counterparts ${\bar{\theta}}^{\dot{B}}$, making the choice consistent with the introduced $Z_6$-grading
$$\theta^A {\bar{\theta}}^{\dot{B}} = - j \; {\bar{\theta}}^{\dot{B}} \theta^A, \; \; \; \;  
{\bar{\theta}}^{\dot{B}}  \theta^A = -j^2 \; \theta^A {\bar{\theta}}^{\dot{B}}.$$

Consider a tri-linear form $\rho^{\alpha}_{ABC}$. We shall call this form $Z_3$-invariant if
we can write:
$$
\rho^{\alpha}_{ABC} \, \theta^A  \theta^B  \theta^C = \frac{1}{3} \, \biggl[ \rho^{\alpha}_{ABC} \,
\theta^A  \theta^B  \theta^C
+ \rho^{\alpha}_{BCA} \, \theta^B  \theta^C  \theta^A +  \rho^{\alpha}_{CAB} \, \theta^C  \theta^A  \theta^B \biggr] =$$
\begin{equation}
= \frac{1}{3} \, \biggl[ \rho^{\alpha}_{ABC} \, \theta^A  \theta^B  \theta^C
+ \rho^{\alpha}_{BCA} \, (j^2 \, \theta^A  \theta^B  \theta^C) +  \rho^{\alpha}_{CAB} \,
j \, (\theta^A  \theta^B  \theta^C) \biggr],
\label{defrhomatrix1}
\end{equation}
by virtue of the commutation relations (\ref{ternary1}).

>From this it follows that we should have
\begin{equation}
\rho^{\alpha}_{ABC} \, \, \theta^A  \theta^B  \theta^C  = \frac{1}{3} \, \biggl[ \rho^{\alpha}_{ABC}
+ j^2 \, \rho^{\alpha}_{BCA} + j \,  \rho^{\alpha}_{CAB} \biggr] \, \theta^A  \theta^B  \theta^C ,
\label{defrhomatrix2}
\end{equation}
from which we get the following properties of the $\rho$-cubic matrices:
\begin{equation}
\rho^{\alpha}_{ABC} = j^2 \, \rho^{\alpha}_{BCA} = j \,  \rho^{\alpha}_{CAB}.
\label{defrhomatrix3}
\end{equation}

Even in this minimal and discrete case, there are covariant and contravariant
indices: the lower and the upper indices display the inverse transformation property. If a given
cyclic permutation is represented by a multiplication by $j$ for the upper indices, the same permutation performed
on the lower indices is represented by multiplication by the inverse, i.e. $j^2$, so that they
compensate each other.

Similar reasoning leads to the definition of the conjugate forms
$ {\bar{\rho}}^{{\dot{\alpha}}}_{{\dot{C}}{\dot{B}}{\dot{A}}}$
satisfying the relations similar to (\ref{defrhomatrix3}) with $j$ replaced be its conjugate, $j^2$:
\begin{equation}
{\bar{\rho}}^{{\dot{\alpha}}}_{{\dot{A}}{\dot{B}}{\dot{C}}} = j \,
{\bar{\rho}}^{{\dot{\alpha}}}_{{\dot{B}}{\dot{C}}{\dot{A}}}
= j^2 \, {\bar{\rho}}^{{\dot{\alpha}}}_{{\dot{C}}{\dot{A}}{\dot{B}}}
\label{defrhomatrix4}
\end{equation}

In the simplest case of two generators, the $j$-skew-invariant forms have
only two independent components:
$$\rho^{1}_{121} = j \, \rho^{1}_{211}
= j^2 \, \rho^{1}_{112},$$
$$\rho^{2}_{212} = j \, \rho^{2}_{122}
= j^2 \, \rho^{2}_{221},$$
and we can set
$$\rho^{1}_{121} = 1, \, \, \rho^{1}_{211} = j^2,  \, \, \rho^{1}_{112} = j,$$
$$\rho^{2}_{212} = 1, \,  \, \rho^{2}_{122} = j^2, \,  \, \rho^{2}_{221} = j.$$

The constitutive cubic relations between the generators of the $Z_3$ graded algebra
can be considered as intrinsic if they are conserved after linear transformations with commuting
(pure number) coefficients, i.e. if they are independent of the choice of the basis.

Let $U^{A'}_A$ denote a non-singular $N \times N$ matrix, transforming the generators
$\theta^A$ into another set of generators, $\theta^{B'} = U^{B'}_B \, \theta^B$.
In principle, the generators of the $Z_2$-graded subalgebra $\xi^{\alpha}$ may or may not
undergo a change of basis. Uniting the two subalgebras in one $Z_6$-graded algebra
suggests that a change of basis should concern all generators at once, both $\xi^{\alpha}$ and $\theta^A$,
This means that under the simultaneous change of basis,
\begin{equation}
\xi^{\alpha} \rightarrow {\tilde{\xi}}^{\beta'} = S^{\beta'}_{\alpha} \; \xi^{\alpha} \; \xi^{\alpha}, \; \; \; \; 
\theta^A \rightarrow {\tilde{\theta}}^{B'} = U^{B'}_A \; \theta^A,
\label{bothxitheta}
\end{equation} 
It seems natural to identify the upper indices $\alpha, \beta$ appearing in the $\rho$-tensors with the
indices appearing in the generators $\xi^{\alpha}$ of the $Z_2$-graded subalgebra. Therefore,
we are looking for the solution of the simultaneous invariance condition for the  $\epsilon_{\alpha \beta}$ and  
$\rho^{\alpha}_{ABC} $ tensors:
\begin{equation}
\epsilon_{\alpha' \beta'} = S^{\alpha}_{\alpha'} S^{\beta}_{\beta'} \; \epsilon_{\alpha \beta},
\; \; \; \; \; S^{{\alpha}'}_{\beta} \, \rho^{{\beta}}_{ABC} = U^{A'}_{A} \, U^{B'}_B \, U^{C'}_C \,
\rho^{{\alpha}'}_{A' B' C'},
\label{covtrans1}
\end{equation}
so that in new basis the numerical values of borth tensors remain the same as before, just like the components
of the Minkowskian space-time metric tensor $g_{\mu \nu} $ remain unchanged under the Lorentz transformations.
Notice that in the last formula above, (\ref{covtrans1}), the matrix $S^{\alpha'}_{\alpha}$ is the inverse matrix
for $S^{\alpha}_{\alpha'}$ appearing in the transformation of the basis $\xi^{\beta}$.

Now, $\rho^{1}_{121} = 1$, and we have two equations corresponding to the choice of values of the index $\alpha'$
equal to $1$ or $2$. For $\alpha' = 1'$ the $\rho$-matrix on the right-hand side is $\rho^{1'}_{A' B' C'}$,
which has only three components,
$$\rho^{1'}_{1' 2' 1'}=1, \, \ \ \, \rho^{1'}_{2' 1' 1'}=j^2, \, \ \ \, \rho^{1'}_{1' 1' 2'}=j, $$
which leads to the following equation:

{\small
\begin{equation}
S^{1'}_{1} = U^{1'}_{1} \, U^{2'}_2 \, U^{1'}_1 + j^2 \, U^{2'}_{1} \, U^{1'}_2 \, U^{1'}_1
+ j \, U^{1'}_{1} \, U^{1'}_2 \, U^{2'}_1 = U^{1'}_{1} \, (U^{2'}_2 \, U^{1'}_1 - U^{2'}_{1} \, U^{1'}_2),
\label{invariant1}
\end{equation} }
because $j^2 + j = - 1$.
For the alternative choice $\alpha' = 2'$ the $\rho$-matrix on the right-hand side is $\rho^{2'}_{A' B' C'}$,
whose three non-vanishing components are
$$\rho^{2'}_{2' 1' 2'}=1, \, \ \ \, \rho^{2'}_{1' 2' 2'}=j^2, \, \ \ \, \rho^{2'}_{2' 2' 1'}=j. $$
The corresponding equation becomes now:
{\small
\begin{equation}
S^{2'}_{1} = U^{2'}_{1} \, U^{1'}_2 \, U^{2'}_1 + j^2 \, U^{1'}_{1} \, U^{2'}_2 \, U^{2'}_1
+ j \, U^{2'}_{1} \, U^{2'}_2 \, U^{1'}_1 = U^{2'}_{1} \, (U^{1'}_2 \, U^{2'}_1 - U^{1'}_{1} \, U^{2'}_2),
\label{invariant2}
\end{equation} }
The two remaining equations are obtained in a similar manner. We choose now the three lower indices
on the left-hand side equal to another independent combination, $(212)$. Then the $\rho$-matrix on the
left hand side must be $\rho^2$ whose component $\rho^2_{212}$ is equal to $1$. This leads to the
following equation when $\alpha' = 1'$:
{\small
\begin{equation}
S^{1'}_{2} = U^{1'}_{2} \, U^{2'}_1 \, U^{1'}_2 + j^2 \, U^{2'}_{2} \, U^{1'}_1 \, U^{1'}_2
+ j \, U^{1'}_{2} \, U^{1'}_1 \, U^{2'}_2 = U^{1'}_{2} \, (U^{1'}_2 \, U^{2'}_1 - U^{1'}_{1} \, U^{2'}_2),
\label{invariant3}
\end{equation}
and the fourth equation corresponding to $\alpha' = 2'$ is:
\begin{equation}
S^{2'}_{2} = U^{2'}_{2} \, U^{1'}_1 \, U^{2'}_2 + j^2 \, U^{1'}_{2} \, U^{2'}_1 \, U^{2'}_2
+ j \, U^{2'}_{2} \, U^{2'}_1 \, U^{1'}_2 = U^{2'}_{2} \, (U^{1'}_1 \, U^{2'}_2 - U^{2'}_{1} \, U^{1'}_2).
\label{invariant4}
\end{equation}}
%{The invariance group of cubic matrices}
%
%The determinant of the $2 \times 2$ complex matrix $U^{A'}_B$ appears everywhere on the right-hand side.
%The obvious solution relating {\it linearly} the matrices $\Lambda^{\alpha'}_{\beta}$ to the matrices
%$U^{A'}_B$ is to impose
%\begin{equation}
%det \, (U^{A'}_B ) = U^{1'}_1 \, U^{2'}_2 - U^{2'}_{1} \, U^{1'}_2 = 1
%\label{detU}
%\end{equation}
%Then we have
%\begin{equation}
%\Lambda^{1'}_{1} = U^{1'}_1, \, \ \, \Lambda^{2'}_{2} = U^{2'}_2, \, \ \
%\Lambda^{1'}_{2} = - U^{1'}_2, \, \ \, \Lambda^{2'}_{1} = - U^{2'}_1,
%\label{lambdaUmatrices}
%\end{equation}
%from which it follows immediately that also $det \, \Lambda = 1$,
%\begin{equation}
%det \, (\Lambda^{\alpha'}_{\beta} ) = \Lambda^{1'}_1 \, \Lambda^{2'}_2 - \Lambda^{2'}_{1} \, \Lambda^{1'}_2 = 1
%\label{detlambda}
%\end{equation}
%Both conditions (\ref{detU}) and (\ref{detlambda}) define the $SL(2, {\bf C})$ group, the covering group
%of the Lorentz group.
The determinant of the $2 \times 2$ complex matrix $U^{A'}_B$ appears everywhere on the right-hand side.
\begin{equation}
S^{2'}_{1} = - U^{2'}_{1} \, [ {\rm det} (U)],
\label{invariantdet}
\end{equation}
The remaining two equations are obtained in a similar manner, resulting in the following:
\begin{equation}
S^{1'}_{2} = - U^{1'}_{2} \, [ {\rm det} (U)], \, \ \ \, \ \ S^{2'}_{2} = U^{2'}_{2} \, [{\rm det} (U)].
\label{invariant34}
\end{equation}
The determinant of the $2 \times 2$ complex matrix $U^{A'}_B$ appears everywhere on the right-hand side.
Taking the determinant of the matrix $S^{{\alpha}'}_{\beta}$ one gets immediately
\begin{equation}
{\rm det} \, (S) = [ {\rm det} \, (U) ]^3.
\label{detLambdaU}
\end{equation}

However, the $U$-matrices on the right-hand side are defined only up to the phase, which due to the
cubic character of the covariance relations
(\ref{invariant1} - \ref{invariant34}),
and they can take on three different values:
$1$, $j$ or $j^2$, i.e. the matrices $j \, U^{A'}_B$ or $j^2 \,  U^{A'}_B$ satisfy the same relations
as the matrices $U^{A'}_B$ defined above.
The determinant of $U$ can take on the values $1, \, j \,$ or $j^2$ if  $det (S) = 1$

Another reason to impose the unitarity condition is as follows. It can be
derived if we require the same behavior for the duals, $\rho_{\beta}^{DEF}$. This extra
condition amounts to the invariance of the anti-symmetric tensor $\epsilon^{AB}$, and
this is possible only if the determinant of $U$-matrices is $1$ (or $j$ or $j^2$, because
only cubic combinations of these matrices appear in the tranformation law for $\rho$-forms.

We have determined the invariance group for the simultaneous change of the basis in our $Z_6$-graded algebra.
However, these transformations based on the $SL(2, {\bf C}$ groups combined with complex representation
of the $Z_3$ cyclic group keep invariant the binary constitutive relations between the $Z_2$-graded generators $\xi^{alpha}$
and the ternary constitutive relations between the $Z_3$-graded generators alone, without mentioning their conjugates
${\bar{\xi}}^{\dot{\alpha}}$ and ${\bar{\theta}}^{\dot{B}}$.

Let us put aside for the moment the conjugate $Z_2$-graded variables, and concentrate our attention on the conjugate
$Z_3$-graded generators ${\bar{\theta}}^{\dot{A}}$ and their commutation relations with $\theta^B$ generators, which
we shall modify as:
\begin{equation}
\theta^A \; {\bar{\theta}}^{\dot{B}} = -j \; {\bar{\theta}}^{\dot{B}} \; \theta^A, \; \; \; \; 
{\bar{\theta}}^{\dot{B}} \; \theta^A = - j^2 \; \theta^A \; {\bar{\theta}}^{\dot{B}}.
\label{commutation2}
\end{equation}

 A similar covariance requirement can be formulated with respect to the set of $2$-forms
mapping the quadratic $\theta^A \, {\bar{\theta}}^{\dot{B}}$ combinations into a four-dimensional linear real space.

It is easy to see, by counting the independent combinations of dotted and undotted indices, that the symmetry (\ref{commutation2}) 
imposed on these expressions reduces their number to four: $(1 {\dot{1}}), \; (1, {\dot{2}}), \;  (2 \; {\dot{1}}), \; (2, {\dot{2}})$,
the conjugate combinations of the type $({\dot{A}} \; B)$ being dependent on the first four because of the imposed symmetry properties.

 Let us define two quadratic forms,  $\pi^{\mu}_{A {\dot{B}}}$ and its conjugate ${\bar{\pi}}^{\mu}_{{\dot{B}} A}$
\begin{equation}
\pi^{\mu}_{A {\dot{B}}} \, \theta^A {\bar{\theta}}^{\dot{B}}
\; \; \; {\rm and} \; \; \;  {\bar{\pi}}^{\mu}_{{\dot{B}} A}
\, {\bar{\theta}}^{\dot{B}} \theta^A.
\label{pisymmetry1}
\end{equation}
The Greek indices $\mu, \nu...$ take on four values, and we shall label them
$0,1,2,3$.

The four tensors $\pi^{\mu}_{A {\dot{B}}}$ and their hermitian conjugates
${\bar{\pi}}^{\mu}_{{\dot{B}} A}$ define a bi-linear mapping from the product of
quark and anti-quark cubic algebras into a linear four-dimensional vector space, whose
structure is not yet defined.

Let us impose the following invariance condition:

\begin{equation}
\pi^{\mu}_{A {\dot{B}}} \, \theta^A {\bar{\theta}}^{\dot{B}} = {\bar{\pi}}^{\mu}_{{\dot{B}} A}
{\bar{\theta}}^{\dot{B}} \theta^A.
\label{invbin}
\end{equation}

It follows immediately from (\ref{commutation2}) that
\begin{equation}
\pi^{\mu}_{A {\dot{B}}} = - j^2 \, {\bar{\pi}}^{\mu}_{{\dot{B}} A}.
\label{pisymmetry2}
\end{equation}
Such matrices are non-hermitian, and they can be realized by the following substitution:
\begin{equation}
\pi^{\mu}_{A {\dot{B}}} = j^2 \, i \, {\sigma}^{\mu}_{A {\dot{B}}}, \, \ \ \,
{\bar{\pi}}^{\mu}_{{\dot{B}} A} = - j \, i \, {\sigma}^{\mu}_{{\dot{B}} A}
\label{pidefinition2}
\end{equation}
where ${\sigma}^{\mu}_{A {\dot{B}}}$
are the unit $2 \time 2$ matrix for $\mu = 0$, and the three hermitian Pauli matrices for $\mu = 1,2,3$.

Again, we want to get the same form of these four matrices in another basis. Knowing
that the lower indices $A$ and ${\dot{B}}$ undergo the transformation with matrices $U^{A'}_B$
and ${\bar{U}}^{{\dot{A}}'}_{\dot{B}}$, we demand that there exist some $4 \times 4$ matrices
$\Lambda^{{\mu}'}_{\nu}$ representing the transformation of lower indices by the matrices
$U$ and ${\bar{U}}$:
 \begin{equation}
\Lambda^{{\mu}'}_{\nu} \, \pi^{\nu}_{A {\dot{B}}} = U^{A'}_A \, {\bar{U}}^{{\dot{B}}'}_{\dot{B}}
 \pi^{{\mu}'}_{A' {\dot{B}}'},
\label{pitransform1}
\end{equation}
It is clear that we can replace the matrices  $\pi^{\nu}_{A {\dot{B}}}$ by the corresponding
matrices $\sigma^{\nu}_{A {\dot{B}}}$,
and this defines the vector ($4 \times 4)$ representation of the Lorentz group.

The first four equations relating the $4 \times 4$ real matrices
 $\Lambda^{{\mu}'}_{\nu}$ with the $2 \times 2$ complex matrices
 $U^{A'}_B$  and
 ${\bar U}^{{\dot{A}}'}_{\dot{B}}$  are as follows:

 $$\Lambda^{0'}_0 - \Lambda^{0'}_3 = U^{1'}_2 \, {\bar U}^{{\dot{1}}'}_{\dot{2}}
+ U^{2'}_2 \, {\bar U}^{{\dot{2}}'}_{\dot{2}} $$
 $$\Lambda^{0'}_0 + \Lambda^{0'}_3 = U^{1'}_1 \, {\bar U}^{{\dot{1}}'}_{\dot{1}}
+ U^{2'}_1 \, {\bar U}^{{\dot{2}}'}_{\dot{1}} $$
 $$\Lambda^{0'}_0 - i \Lambda^{0'}_2 = U^{1'}_1 \, {\bar U}^{{\dot{1}}'}_{\dot{2}}
+ U^{2'}_1 \, {\bar U}^{{\dot{2}}'}_{\dot{2}} $$
 $$\Lambda^{0'}_0 + i \Lambda^{0'}_2 = U^{1'}_2 \, {\bar U}^{{\dot{1}}'}_{\dot{1}}
+ U^{2'}_2 \, {\bar U}^{{\dot{2}}'}_{\dot{1}} $$

The next four equations relating the $4 \times 4$ real matrices
 $\Lambda^{{\mu}'}_{\nu}$ with the $2 \times 2$ complex matrices
 $U^{A'}_B$  and
${\bar U}^{{\dot{A}}'}_{\dot{B}}$ are as follows:

 $$\Lambda^{1'}_0 - \Lambda^{1'}_3 = U^{1'}_2 \, {\bar U}^{{\dot{2}}'}_{\dot{2}}
+ U^{2'}_2 \, {\bar U}^{{\dot{1}}'}_{\dot{2}} $$
 $$\Lambda^{1'}_0 + \Lambda^{1'}_3 = U^{1'}_1 \, {\bar U}^{{\dot{2}}'}_{\dot{1}}
+ U^{2'}_1 \, {\bar U}^{{\dot{1}}'}_{\dot{1}} $$
 $$\Lambda^{1'}_1 - i \Lambda^{1'}_2 = U^{1'}_1 \, {\bar U}^{{\dot{2}}'}_{\dot{2}}
+ U^{2'}_1 \, {\bar U}^{{\dot{1}}'}_{\dot{2}} $$
 $$\Lambda^{1'}_1 + i \Lambda^{1'}_2 = U^{1'}_2 \, {\bar U}^{{\dot{2}}'}_{\dot{1}}
+ U^{2'}_2 \, {\bar U}^{{\dot{1}}'}_{\dot{1}} $$
We skip the next two groups of four equations corresponding to the "spatial" indices $2$ and $3$,
reproducing the same scheme as the last four equations with the space index equal to $1$.

It can be checked that now ${\rm det} \;
(\Lambda) = \left[ {\rm det} U \right]^2 \, \left[ {\rm det} {\bar{U}} \right]^2.$

The group of transformations thus defined is $SL(2, {\bf C})$, which is the covering group of the Lorentz group.

%{The metric tensor $g_{\mu \nu}$}

With the invariant ``spinorial metric" in two complex dimensions, $\varepsilon^{AB}$
and $\varepsilon^{{\dot{A}}{\dot{B}}}$ such that $\varepsilon^{12} = - \varepsilon^{21} = 1$
and $\varepsilon^{{\dot{1}}{\dot{2}}} = - \varepsilon^{{\dot{2}}{\dot{1}}}$, we can define
the contravariant components $\pi^{\nu \, \, A {\dot{B}}}$. It is easy to show that the
Minkowskian space-time metric, invariant under the Lorentz transformations, can be defined as
\begin{equation}
g^{\mu \nu} = \frac{1}{2} \biggl[ \pi^{\mu}_{A {\dot{B}}} \, \pi^{\nu \, \, A {\dot{B}}} \biggr]
= diag (+,-,-,-)
\label{Mmetric}
\end{equation}
Together with the anti-commuting spinors ${\psi}^{\alpha}$ the four real coefficients defining
a Lorentz vector, $x_{\mu} \, {\pi}^{\mu}_{A {\dot{B}}}$, can generate now the supersymmetry
via standard definitions of super-derivations.
Let us then choose the matrices $\Lambda^{\alpha'}_{\beta}$ to be the usual spinor representation of
the $SL(2, {\bf C})$ group, while the matrices $U^{A'}_{B}$ will be defined as follows:
\begin{equation}
U^{1'}_{1} = j \Lambda^{1'}_1,   \; \; \; U^{1'}_{2} = - j \Lambda^{1'}_2, \; \; \;
U^{2'}_{1} = - j  \Lambda^{2'}_1, \; \; \;  U^{2'}_{2} = j \Lambda^{2'}_2,
\label{Umatrices}
\end{equation}
the determinant of $U$ being equal to $j^2$.
Obviously, the same reasoning leads to the conjugate cubic representation of $SL(2, {\bf C})$ if we require
the covariance of the conjugate tensor
$${\bar{\rho}}^{\dot{\beta}}_{{\dot{D}}{\dot{E}}{\dot{F}}} = j \,
{\bar{\rho}}^{\dot{\beta}}_{{\dot{E}}{\dot{F}} {\dot{D}}}
= j^2 \, {\bar{\rho}}^{\dot{\beta}}_{{\dot{F}} {\dot{D}} {\dot{E}}},$$
by imposing the equation similar to (\ref{covtrans1})
\begin{equation}
\Lambda^{{\dot{\alpha}}'}_{{\dot{\beta}}} \, {\bar{\rho}}^{{\dot{\beta}}}_{{\dot{A}}{\dot{B}}{\dot{C}}} =
{\bar{\rho}}^{{\dot{\alpha}}'}_{{\dot{A}}' {\dot{B}}' {\dot{C}}'} {\bar{U}}^{{\dot{A}}'}_{{\dot{A}}} \,
{\bar{U}}^{{\dot{B}}'}_{{\dot{B}}} \, {\bar{U}}^{{\dot{C}}'}_{{\dot{C}}} .
\label{covtrans2}
\end{equation}
The matrix $\bar{U}$ is the complex conjugate of the matrix $U$, with determinant equal to $j$.

%Moreover, the two-component entities obtained as images of cubic combinations of quarks,
%$\psi^{\alpha} = \rho^{\alpha}_{ABC} \theta^A \theta^B \theta^C$
% and ${\bar{\psi}}^{\dot{\beta}} = {\bar{\rho}}^{{\dot{\beta}}}_{{\dot{D}}{\dot{E}}{\dot{F}}}
%{\bar{\theta}}^{\dot{D}} {\bar{\theta}}^{\dot{E}} {\bar{\theta}}^{\dot{F}} $ should anti-commute,
%because their arguments do so, by virtue of (\ref{commutation2}):
%$$ (\theta^A \theta^B \theta^C) ({\bar{\theta}}^{\dot{D}} {\bar{\theta}}^{\dot{E}} {\bar{\theta}}^{\dot{F}} )
%= - ({\bar{\theta}}^{\dot{D}} {\bar{\theta}}^{\dot{E}} {\bar{\theta}}^{\dot{F}})(\theta^A \theta^B \theta^C)$$

In conclusion, we have found the way to derive the covering group of the Lorentz group acting on spinors via
the usual spinorial representation, and on vectors via the $4 \times 4$ real matrices.. Spinors are obtained as a homomorphic image of
tri-linear combinations of three $Z_3$-graded generators $\theta^A$ (or their conjugates ${\bar{\theta}}^{\dot{B}}$). The 
$Z_3$-graded generators transform with
matrices $U$ (or ${\bar{U}}$ for the conjugates), but these matrices are not unitary:
their determinants are equal to $j^2$ or $j$, respectively.

In our forthcoming paper we shall investigate similar $Z_6$-graded generalization extended to the differential forms $d \xi^{\alpha}$ and $d \theta^B$.

\vskip 0.3cm
%This does not preclude that the components of $\eta^{\delta'} \eta^{\alpha'}$ are the same as the components of
\noindent
{\small {\bf Acknowledgement}
We express our thanks to Michel Dubois-Violette for enlightening discussions and very useful remarks and suggestions.
V. Abramov and O. Liivapuu gratefully appreciate the financial support of Estonian Science Foundation under the Research Grant No. ETF9328. 
They also express their gratitude to the Estonian Ministry of Education and Research for financial support by institutional research funding IUT20-57.}

\end{document}